\documentclass[11pt]{article}

\usepackage{amssymb}

\def\emline#1#2#3#4#5#6{%
       \put(#1,#2){\special{em:moveto}}%
       \put(#4,#5){\special{em:lineto}}}
\def\newpic#1{}

\newtheorem{predfn}{\bf{Definition}}

\newenvironment{dfn}{\begin{predfn}{\hspace{
               -.5em}{.}}}{\end{predfn}}
\newtheorem{preex}{\bf{Example}}

\newtheorem{prelem}{\bf{Lemma}}

\newenvironment{lem}{\begin{prelem}{\hspace{
               -.5em}{.}}}{\end{prelem}}
\newtheorem{prethm}{\bf{Theorem}}

\newenvironment{thm}{\begin{prethm}{\hspace{
               -.5em}{.}}}{\end{prethm}}
\newtheorem{preathm}{\bf{Theorem}}

\newenvironment{athm}{\begin{preathm}{\hspace{
               -.5em}{.\hspace{.2em}}}}{\end{preathm}}
\newtheorem{prepro}{\bf{Proposition}}

\newenvironment{pro}{\begin{prepro}{\hspace{
               -.5em}{.}}}{\end{prepro}}
\newtheorem{precor}{\bf{Corollary}}

\newtheorem{preconj}{\bf{Conjecture}}

\newenvironment{conj}{\begin{preconj}{\hspace{
               -.5em}{.}}}{\end{preconj}}
\newtheorem{preconjj}{\bf{Conjecture}}

\newenvironment{conjj}{\begin{preconjj}{\hspace{
               -.9em}{.\hspace{.2em}}}}{\end{preconjj}}
\newtheorem{prerem}{\bf{Remark}}

\newtheorem{preproof}{{\bf Proof.}}

\newenvironment{proof}[1]{\begin{preproof}{\rm
               #1}\hfill{$\blacksquare$}}{\end{preproof}}

\newcommand{\m}{\mbox{\rm m}}
\newcommand{\utlc}{U$2$LC}
\newcommand{\uthreelc}{U$3$LC}
\newcommand{\uklc}{U$k$LC}
\newcommand{\uflc}{U$f$LC}
\newcommand{\mnum}{{\rm m}--number}
\newcommand{\ch}{\Chi_\ell}
\newcommand{\ech}{\Chi'_\ell}

\newcommand{\lc}{list critical}
\newcommand{\klc}{$k$--list critical}
\newcommand{\threelc}{$3$--list critical}
\newcommand{\deffnt}{\sf}
\newcommand{\reffnt}{\rm}
\newcommand{\N}{\Bbb{N}}
\def\Chi{\lower-.3ex\hbox{$\chi$}}
\def\le{\leqslant}
\def\ge{\geqslant}

\title{Some Concepts in List Coloring} %
\author{
Ch.~Eslahchi{$^\ast$}, M.~Ghebleh{$^\dag$}, and H.~Hajiabolhassan{$^\ast$}
\vspace{1eX}\\ %
\normalsize Institute for Studies in %
Theoretical Physics and Mathematics (IPM)\vspace{2.3eX}\\ %
$^\ast$\normalsize Department of Mathematics\\
\normalsize Shahid Beheshti University\\
\normalsize Evin, Tehran, Iran\vspace{2eX}\\
$^\dag$\normalsize Department of Mathematical Sciences\\
\normalsize Sharif University of Technology\\
\normalsize P.\,O.~Box~11365--9415, Tehran, Iran\\
} %
\date{}

\begin{document}
\addtolength{\baselineskip}{1mm}
\maketitle

\begin{abstract}
In this paper uniquely list colorable graphs are studied. A graph
$G$ is called to be uniquely $k$--list colorable if it admits
a $k$--list assignment from which $G$ has a unique list coloring.
The minimum $k$ for which $G$ is not uniquely $k$--list colorable
is called the \mnum\ of $G$.
We show that every triangle--free uniquely colorable graph with
chromatic number $k+1$, is uniquely $k$--list colorable. A bound
for the \mnum\ of graphs is given, and using this bound it is shown
that every planar graph has \mnum\ at most $4$. Also we introduce
list criticality in graphs and characterize all \threelc\ graphs.
It is conjectured that every $\ech$--critical graph is
$\Chi'$--critical and the equivalence of this conjecture to the
well known list coloring conjecture is shown.
\end{abstract}

\section{Introduction}

We consider finite, undirected simple graphs. For necessary
definitions and notations we refer the reader to standard texts
such as \cite{west}.

By a {\deffnt $k$--list assignment} $L$ to a graph $G$ we mean a
map which assigns to each vertex $v$ of $G$ a set $L(v)$ of size
$k$. A {\deffnt list coloring} for $G$ from $L$, or an {\deffnt
$L$--coloring} for short, is a proper coloring $c$, in which for
each vertex $v$, $c(v)$ is chosen from $L(v)$. A graph $G$ is
called {\deffnt $k$--choosable} if it has a list coloring from
any $k$--list assignment to it. The minimum number $k$
for which $G$ is $k$--choosable is called the {\deffnt list
chromatic number} of $G$ and is denoted by $\ch(G)$. In the
following theorem all $2$--choosable graphs are characterized.
Before we state the theorem it should be noted that the core of
a graph is a subgraph which is obtained by repeatedly deleting
a vertex of degree 1, until no vertex of degree 1 remains.

\begin{athm}\label{2choose}{\reffnt\cite{erdos}}
A connected graph is $2$--choosable, if and only if its core is
either a single vertex, an even cycle, or $\theta_{2,2,2r}$, for
some $r\ge 1$.
\end{athm}

A graph $G$ is called {\deffnt uniquely $k$--list colorable}, or
{\deffnt\uklc} for short, if it admits a $k$--list assignment $L$
such that $G$ has a unique $L$--coloring. This
concept was introduced by Dinitz and Martin~\cite{dinmar} and
independently by Mahdian and Mahmoodian~(\cite{mahmah} and
\cite{tabriz}). A characterization of uniquely $2$--list colorable
graphs follows.

\begin{athm}\label{u2lc}{\reffnt\cite{mahmah}}
A graph $G$ is not \utlc\ if and only if each of its blocks is
either a cycle, a complete graph, or a complete bipartite graph.
\end{athm}

It is easy to see that for each graph $G$ there exists a number
$k$ such that $G$ is not \uklc. The minimum $k$ with this property
is called the {\deffnt\mnum} of $G$ and is denoted by $\m(G)$. It
is shown in~\cite{tabriz} that every planar graph has \mnum\ at most
$5$, and it is asked about the existence of planar graphs with \mnum\
equal to $5$. We study uniquely list colorable graphs in
Section~\ref{uklcsec}, where we prove that every triangle--free
uniquely $(k+1)$--colorable graph is uniquely $k$--list colorable.
We also show that every planar graph has \mnum\ at most $4$, so
the answer to that question in~\cite{tabriz} is negative.

In Section~\ref{criticalsec} we introduce list critical graphs and
characterize $3$--list critical graphs. Finally we pose a conjecture
about list critical graphs which is shown to be equivalent to the
list coloring conjecture.

\section{Uniquely list colorable graphs\label{uklcsec}}

In~\cite{ghbmah} one can find several examples of \uklc\ graphs,
for some arbitrary positive integer $k$. In the following
lemma we also introduce a class of \uklc\ graphs. In this way we
relate uniquely list colorable graphs to uniquely colorable
graphs.

\begin{lem}
\label{uvc_uklc} Let $G$ be a uniquely colorable graph with chromatic
number $k+1$, and $c$ be its unique $(k+1)$--coloring with color
classes $C_1,\ldots, C_{k+1}$. If for each $i\le k+1$, $|C_i|\ge i-1$,
then $G$ is a uniquely $k$--list colorable graph.
\end{lem}

\begin{proof}{
We proceed by induction on $k$ and prove that there exists a $k$--list
assignment to such graph $G$ using exactly $k+1$ colors, which induces
a unique list coloring. For $k=1$ the result obviously holds.
Let $G$ be a uniquely $(k+1)$--colorable graph as in the statement and
$k\ge 2$. By induction $G\setminus C_{k+1}$ admits a $(k-1)$--list
assignment $L'$ which induces a unique list coloring and uses colors
$1,\ldots,k$. For each $v\in V(G)\setminus C_{k+1}$, assign the list
$L(v)=L'(v)\cup\{k+1\}$ to $v$, and since $|C_{k+1}|\ge k$, it is
possible to assign some lists to $C_{k+1}$ such that
$\bigcap_{v\in C_{k+1}} L(v)=\{k+1\}$. Now it is easy to see that $L$
is the desired list assignment.
}\end{proof}

It is shown in {\reffnt\cite{wangartzy}} that for every $k\ge 3$,
in a triangle--free uniquely $k$--colorable graph, each color
class has at least $k+1$ vertices. Using this result, we obtain
the following theorem.

\begin{thm}
Every triangle--free uniquely $(k+1)$--colorable graph is uniquely
$k$--list colorable.
\end{thm}

On the other hand in~{\reffnt\cite{largegirth}} it is shown that
for each $k\ge 2$, there exists a uniquely $k$--colorable graph
with arbitrary large girth. So by theorem above, for each $k$, there
exists a \uklc\ graph with arbitrary large girth.

We need here a definition which is a generalization of the concept
of a \uklc\ graph.

\begin{dfn}
Let $G$ be a graph and $f$ a be function from $V(G)$ to $\N$. An
{\deffnt $f$--list assignment} $L$ to $G$ is a list assignment in
which $|L(v)|=f(v)$ for each vertex $v$. The graph $G$
is called to be {\deffnt uniquely $f$--list colorable}, or
{\deffnt\uflc} for short, if there exists an $f$--list assignment
$L$ for it such that $G$ has a unique $L$--coloring.
\end{dfn}

By definition above, if $G$ is a \uflc\ graph, where $f(v)=k$ for
each vertex $v$ of $G$, then $G$ in fact is a \uklc\ graph.
To prove the next theorem, we need a relation which is proved in
Truszczy\'nski~{\reffnt\cite{trusz}} and states that if $G$ is a
uniquely $k$--colorable graph, then
$e(G)\ge (k-1)n(G)-{k\choose 2}.$ %

\begin{thm}
\label{sigmafv} If $G$ is a \uflc\ graph,
then $$\sum_{v\in V(G)}f(v)\le n(G)+e(G).$$
\end{thm}

\begin{proof}{
Suppose that $L$ is an $f$--list assignment to $G$
using colors $1,2,\ldots,t$, such that $G$ has a unique
$L$--coloring. We construct a uniquely $t$--colorable graph
$G^\ast$ as follows. Let $V(G)=\{v_1,\ldots,v_n\}$ and $K_t$ is a
complete graph on the vertex set $\{w_1,\ldots,w_t\}$. Now for
$G^\ast$ consider the union of $G$ and $K_t$ and add edges
$v_iw_j$ where $1\le i\le n$, $1\le j\le t$, and $j\not\in
L(v_i)$.

Consider a $t$--coloring $c$ of $G^\ast$. Without loss of
generality we can assume that $c(w_i)=i$ for each $1\le i\le t$.
Since $G$ has a unique $L$--coloring, by construction of $G^\ast$,
$c$ is the only $t$--coloring of $G^\ast$. So $G^\ast$ is a
uniquely $t$--colorable graph. On the other hand $G^\ast$ has
$n(G)+t$ vertices, and $e(G)+{t\choose 2}+\sum_{v\in
V(G)}(t-f(v))$ edges. Therefore as mentioned above, we have
$$e(G)+{t\choose 2}+\sum_{v\in V(G)}(t-f(v))\ge
(n(G)+t)(t-1)-{t\choose 2}$$ and after
simplification we obtain the result.
}\end{proof}

A natural question which arises here is that whether or not
equality holds in Theorem~\ref{sigmafv}? In the following
proposition we give a positive answer to this question.

\begin{pro}
For every graph $G$, there exists $f:V(G)\to\N$ such that $G$ is
\uflc\ and $\sum_{v\in V(G)}f(v)=n(G)+e(G)$.
\end{pro}

\begin{proof}{
We proceed by induction on the number of vertices of $G$. For
$n(G)=1$ the statement is obvious. Consider a graph $G$ with
$n(G)\ge 2$ and a vertex $v$ of $G$. By induction there exists
$f':V(G\setminus v)\to \N$ and an $f'$--list assignment $L'$ to
$G\setminus v$ such that $G\setminus v$ has a
unique $L'$--coloring, and $\sum_{w\in V(G\setminus
v)}f(w)=n(G\setminus v)+e(G\setminus v)$. Consider a color $a$
which is not used by $L'$, and define a list assignment $L$ to
$G$ as follows. %
$$L(w)=\left\{\begin{array}{ll} %
a&{\rm for\ }w=v\\ %
L'(w)\cup\{a\}&{\rm for\ }w\in N(v)\\ %
L'(w)&{\rm otherwise.}\\ %
\end{array}\right.$$ %
It is easy to verify that $G$ has a unique $L$--coloring and that
we have $\sum_{v\in V(G)}|L(v)|=n(G)+e(G)$.
}\end{proof}

Although the proposition above shows that in Theorem~\ref{sigmafv}
equality may hold, but it seems that if $f(v)=k$ for each vertex
$v$, equality does not hold and we have $e(G)>(k-1)n(G)$.

By definition every graph $G$ for $k=\m(G)-1$ is \uklc\. So by
Theorem~\ref{sigmafv} we have the following.

\begin{thm}
\label{bound} For a graph $G$ let $\overline{d}(G)$ denote the
average degree of $G$, i.e. $\overline{d}(G)=2e(G)/n(G)$. Then
$$\m(G)\le \lfloor{\overline{d}(G)\over 2}\rfloor+2.$$
\end{thm}

For example suppose that $G$ is a bipartite graph. We have
$\overline{d}(G)\le n(G)/2$ so Theorem~\ref{bound} implies
$\m(G)\le\lfloor n(G)/4 +2\rfloor$. This bound can be improved
to a logarithmic bound as we will show in Theorem~\ref{logbnd},
but first we need a lemma.

Let $L$ be a $k$--list assignment to a graph $G$ such that $G$ has
a unique $L$--coloring $c$. For each vertex $v$ of $G$, all the
elements of $L(v)\setminus\{c(v)\}$ must appear in $N(v)$, so if
we denote by $c_N(v)$ the set of colors appearing in $N(v)$, then
$|c_N(v)|\ge k-1$. In the following lemma we state a stronger result.

\begin{lem}
\label{kcolor} Suppose that $G$ is a \uklc\ graph, and $L$ is a
$k$--list assignment to $G$ such that $G$ has a unique $L$--coloring
$c$ with color classes $C_1,\ldots,C_t$ such that $c(C_i)=\{i\}$.
There exist at least $k-1$ classes containing a vertex~$v$ with
$|c_N(v)|\ge k$.
\end{lem}

\begin{proof}{
Without loss of generality suppose that for $\ell\ge k-1$,
$C_\ell$ contains no vertex $v$ with $|c_N(v)|\ge k$. Assume that
$u\in C_{k-1}$, $i=c(v_0)$, and $j\in
L(v_0)\setminus\{1,\ldots,k-1\}$. Suppose that $G_{ij}$ is the
subgraph of $G$ induced on $C_i\cup C_j$. Since for each vertex
$v$ of the component of $G_{ij}$ containing $v_0$ we have
$|c_N(v)|=k-1$, it is implied that $i,j\in L(v)$. So we can
interchange the colors $i$ and $j$ in this component to obtain a
new $L$--coloring for $G$. This contradiction completes the proof.
}\end{proof}

It is shown in~\cite{erdos} that every non--$k$--choosable bipartite
graph has more than $2^{k-1}$ vertices. So by applying
Lemma~\ref{kcolor}, we deduce the following theorem.

\begin{thm}
\label{logbnd}
Let $G$ be a bipartite graph. Then $\m(G)\le 2+\log_2n(G).$
\end{thm}

\begin{proof}{
Suppose that $L$ is a $k$--list assignment to $G$
such that $G$ has a unique $L$--coloring $c$. By
Lemma~\ref{kcolor}, $G$ has a vertex $v_0$, such that there are at
least $k$ colors appeared at $N(v_0)$ in $c$. Let $G'$ be the
graph obtained from $G$, by duplicating $v_0$, i.e. adding a new
vertex $w$ to $G$ and joining it to $N(v_0)$. Now assign to $w$ a
list containing $k$ of the colors appeared at $N(v_0)$ in $c$, and
the list $L(v)$ to each other vertex $v$ of $G'$. It is clear that
$G'$ is a bipartite graph and it has no coloring from these lists,
so it is not $k$--choosable. Hence $n(G')>2^{k-1}$ vertices. This
implies that $n(G)\ge 2^{k-1}$, and so $k\le 1+\log_2n(G)$. Now we
obtain the desired relation by setting $k=\m(G)-1$.
}\end{proof}

In the remainder of this section we state some consequences of
Theorems~\ref{sigmafv} and \ref{bound}.

It is well known that a
planar graph with $n$~vertices has at most $3n-6$ edges. So the
following theorem is an immediate consequence of
Theorem~\ref{bound}.

\begin{thm}
\label{plnrbbnd} For every planar graph $G$ we have $\m(G)\le 4$.
\end{thm}

\begin{figure} %
\centering %
\special{em:linewidth 0.4pt} \unitlength .70mm %
\linethickness{0.4pt} %
\begin{picture}(82.92,57.92) %
\put(1.67,1.67){\circle*{2.50}} %
\put(81.67,1.67){\circle*{2.50}} %
\put(41.67,16.67){\circle*{2.50}} %
\put(41.67,26.67){\circle*{2.50}} %
\put(41.67,36.67){\circle*{2.50}} %
\put(41.67,46.67){\circle*{2.50}} %
\put(41.67,56.67){\circle*{2.50}} %
\emline{41.67}{56.67}{1}{41.67}{16.67}{2} %
\emline{41.67}{16.67}{3}{1.67}{1.67}{4} %
\emline{1.67}{1.67}{5}{81.67}{1.67}{6} %
\emline{81.67}{1.67}{7}{41.67}{16.67}{8} %
\emline{1.67}{1.67}{9}{41.67}{26.67}{10} %
\emline{41.67}{26.67}{11}{81.67}{1.67}{12} %
\emline{81.67}{1.67}{13}{41.67}{36.67}{14} %
\emline{41.67}{36.67}{15}{1.67}{1.67}{16} %
\emline{1.67}{1.67}{17}{41.67}{46.67}{18} %
\emline{41.67}{46.67}{19}{81.00}{1.67}{20} %
\emline{81.00}{1.67}{21}{41.67}{56.67}{22} %
\emline{41.67}{56.67}{23}{1.67}{1.67}{24} %
\end{picture}
\caption{A uniquely $3$--list colorable planar graph} %
\label{u3lc}
\end{figure}
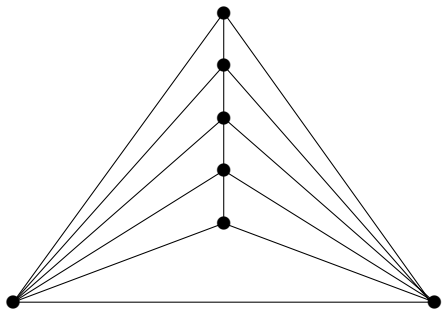

By Lemma~\ref{uvc_uklc} the planar graph shown in Figure~\ref{u3lc}
is a \uthreelc\ graph for which the inequalities in
Theorem~\ref{bound} and Theorem~\ref{plnrbbnd} turn to be equalities.

Furthermore we know that a triangle--free planar graph $G$, has at
most $2n(G)-4$ edges. So Theorem~\ref{bound} implies that each
triangle--free planar graph $G$ has \mnum\ at most $3$. In the
following proposition a stronger result is obtained.

\begin{pro}
\label{8face} If a plane graph has at most $7$ triangular faces,
then $\m(G)\le 3$.
\end{pro}

\begin{proof}{
Consider a \uthreelc\ plane graph $G$ with $n$ vertices, $e$
edges, $f$ faces, and $t$ triangular faces. We have $2e\ge
4(f-t)+3t=4f-t$, and by Euler formula $f=2-n+e$, so
$t\ge 8-4n+2e$. On the other hand Theorem~\ref{sigmafv} implies
that $e\ge 2n$. So $t\ge 8$, as desired.
}\end{proof}

The following conjecture is about the structure of \uthreelc\
planar graphs which is motivated by the proposition above.

\begin{conj}
Every \uthreelc\ planar graph has $K_4$ as a subgraph.
\end{conj}

For another application of Theorem~\ref{sigmafv}, we study line
and total versions of uniquely list coloring.

A graph $G$ is called to be {\deffnt uniquely $k$--list edge
colorable}, if $L(G)$ is a uniquely $k$--list colorable graph. The
{\deffnt edge \mnum} of $G$ is defined to be $\m(L(G))$, and is
denoted by $\m'(G)$. It is straightforward to see that for each
graph $G$, $\overline{d}(L(G))\le\Delta(L(G))\le 2\Delta(G)-2$.
So using Theorem~\ref{bound} we deduce the following.

\begin{thm}
\label{edgemnum}
For every graph $G$, we have\ \,$\m'(G)\le\Delta(G)+1$ and if
\ \,$\m'(G)=\Delta(G)+1$ then $G$ is a regular graph.
\end{thm}

Note that in Theorem~\ref{edgemnum} it is shown that if $G$ is not a
regular graph, then $\m(G)\le\Delta(G)$. So in this case
$\m(G)\le\Chi(G)$.

\section{List critical graphs\label{criticalsec}}

In this section we introduce a concept of \lc\ graphs and we
state some results concerning it.

\begin{dfn}
\label{lcrit} A graph $G$ is called $\ch$--critical if for each
proper subgraph $H$ of it we have $\ch(H)<\ch(G)$.
\end{dfn}

We sometimes refer to a $\ch$--critical graph $G$ as a $k$--list
critical graph, where $k=\ch(G)$. It can easily be verified that
the only connected $2$--list critical graph is $K_2$, odd cycles
are $3$--list critical, and the complete graph $K_k$ is $k$--list
critical.

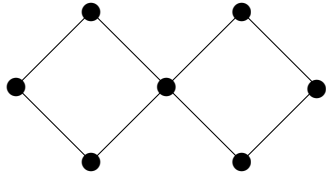
\begin{figure}
\centering %
\special{em:linewidth 0.4pt} %
\unitlength 1.00mm %
\linethickness{0.4pt} %
\begin{picture}(43.25,22.92) %
\put(2.00,11.67){\circle*{2.50}} %
\put(22.00,11.67){\circle*{2.50}} %
\put(42.00,11.34){\circle*{2.50}} %
\put(32.00,21.67){\circle*{2.50}} %
\put(32.00,1.67){\circle*{2.50}} %
\put(12.00,1.67){\circle*{2.50}} %
\put(12.00,21.67){\circle*{2.50}} %
\emline{2.00}{11.67}{1}{12.00}{21.67}{2} %
\emline{12.00}{21.67}{3}{32.00}{1.67}{4} %
\emline{32.00}{1.67}{5}{42.00}{11.67}{6} %
\emline{42.00}{11.67}{7}{32.00}{21.67}{8} %
\emline{32.00}{21.67}{9}{12.00}{1.67}{10} %
\emline{12.00}{1.67}{11}{2.00}{11.67}{12} %
\end{picture}
\caption{A non--$2$--connected $3$--list critical graph} %
\label{tie}
\end{figure}

Obviously every graph $G$ contains a $\ch$--critical subgraph $H$
such that $\ch(H)=\ch(G)$, and by an argument similar to critical
graphs, $\delta(G)\ge \ch(G)-1$. On the other hand there
exists some differences between critical graphs and list critical
graphs. For example it is well known that every critical graph is
$2$--connected. In Figure~\ref{tie} we have given an example of a
$3$--list critical graph which is not $2$--connected.

In the next theorem \threelc\ graphs are characterized.

\begin{thm}
A graph is \threelc\ if and only if it is either
an odd cycle, %
two even cycles with a path joined them, %
$\theta_{r,s,t}$ where $r$, $s$, $t$ have the same parity, and
at most one of them is~$2$, or %
$\theta_{2,2,2,2r}$ where $r\ge 1$. %
\end{thm}

\begin{proof}{
By use of Theorem~\ref{2choose}, it is easy to see that all the
graphs listed in the statement are \threelc.

For the converse suppose that $G$ is a \threelc\ graph.
If $G$ is $2$--connected, by a theorem of Whitney~\cite{whitney}
$G$ has an ear decomposition $K_2\cup P^1\cup\ldots\cup P^q$.
If $q\ge 4$, deleting an edge of $P^q$, yields a non--$2$--choosable
graph, which contradicts the \threelc ity of $G$. So $q\le 3$ and we
consider the following three cases.

\begin{itemize} %
\item If $q=1$, $G$ is a cycle, and so it is an odd cycle.
\item If $q=2$, $G=\theta_{r,s,t}$. In this case by deleting each
edge of $G$, we obtain a graph whose core is a cycle, and since
this cycle must be even, the numbers $r$, $s$, and $t$ have the
same parity. Now if at least two of $r$, $s$, and $t$ are equal to
$2$, we have $\ch(G)=2$, a contradiction.
\item The last case is $q=3$. By deleting each edge of $G$ we
obtain a graph whose core is a $\theta_{r,s,t}$, and since this
graph must be $2$--choosable, we have $r=s=2$ and $t$ is an even
number. Now by case analysis, it is easy to see that
$G=\theta_{2,2,2,2\ell}$.
\end{itemize}
On the other hand if $G$ is not $2$--connected, we consider two
end--blocks $B_1$ and $B_2$ of $G$. Since $\delta(G)\ge 2$ each of
$B_1$ and $B_2$ has a cycle. So $G$ has a subgraph $H$ which is
composed of two edge--disjoint cycles joined to each other by a
path (possibly of length zero). We know that $\ch(H)=3$, and so by
$\ch$--criticality of $G$, $G$ has no edge outside $H$, i.e.
$G=H$. Hence $G$ satisfies the statement.
}\end{proof}

Suppose that $G$ is a $k$--list critical graph, and $L$ is a
$k$--list assignment to $G$. Consider a
vertex $v$ in $G$ and a color $a\in L(v)$. Assign to each vertex
$u$ in $G\setminus v$ the list $L(u)\setminus\{a\}$. Since
$G\setminus v$ is $(k-1)$--choosable, it has a coloring from
the assigned lists, and one can extend this coloring to an
$L$--coloring of $G$ by assigning the color $a$ to $v$. So there
exists an $L$--coloring for $G$ in which $v$ takes $a$.

As mentioned in the previous paragraph, every \klc\ graph has at
least $k$ colorings from each $k$--list assignment
so every \klc\ graph has \mnum\ at most $k$.

A graph $G$ is called to be {\deffnt edge $k$--choosable}, if the
graph $L(G)$ is $k$--choosable, and the {\deffnt list chromatic
index} of $G$ written $\ech(G)$ is defined to be $\ch(L(G))$. As
in the case of defining $\Chi'$--critical graphs, one can define a
$\ech$--critical graph $G$ to be a graph in which for each proper
subgraph $H$, $\ech(H)<\ech(G)$. We recall here the well known
List Coloring Conjecture (LCC), which first appeared in print
in~\cite{LCCref}.

\begin{conjj}{\reffnt\cite{LCCref}}
Every graph $G$ satisfies $\ech(G)=\Chi'(G)$.
\end{conjj}

Suppose that $G$ is a counterexample to the LCC with minimum
number of edges. So for each edge $u v$ of $G$ we have
$\ech(G\setminus uv)=\Chi'(G\setminus uv)$, and since
$\Chi'(G\setminus uv)\le\Chi'(G)<\ech(G)$, we conclude that
$\ech(G\setminus uv)=\ech(G)-1$. This means that $G$ is a
$\ech$--critical graph and therefore $\ech$--critical graphs
may be useful to attack the LCC.

In the study of $\ech$--critical graphs we have lead to the
following conjecture.

\begin{conj}
\label{CLCC} Every $\ech$--critical graph is $\Chi'$--critical.
\end{conj}

\begin{pro}
The conjecture above is equivalent with the LCC,
while its converse is implied by the LCC.
\end{pro}

\begin{proof}{
It is straight forward to check that the list coloring conjecture
implies Conjecture~\ref{CLCC} and its converse. On the other hand
suppose that Conjecture~\ref{CLCC} is true, and $G$ is a
counterexample to the list coloring conjecture with minimum number
of edges. As mentioned above $G$ is $\ech$--critical, and by
Conjecture~\ref{CLCC}, it is $\Chi'$--critical. By removing an
arbitrary edge $uv$ from $G$ we obtain a graph for which the list
coloring conjecture holds. So $\ech(G\setminus
uv)=\Chi'(G\setminus uv)$, and this means that
$\ech(G)-1=\Chi'(G)-1$, a contradiction.
}\end{proof}

In~\cite{galvin} it is proved that every bipartite multigraph
fulfills the LCC. On the other hand we know that the only
bipartite $\Chi'$--critical graphs are stars. So the following
theorem is implied by a similar argument as in the previous
paragraph, the only bipartite $\ech$--critical graphs
are stars.

\section*{Acknowledgements}
The authors wish to thank Professor E.\,S.~Mahmoodian and
R.~Tusserkani for their helpful comments. They are indebted to the
Institute for Studies in Theoretical Physics and Mathematics
(IPM), Tehran, Iran for their support.
\newcommand{\noopsort}[1]{} \newcommand{\printfirst}[2]{#1}
  \newcommand{\singleletter}[1]{#1} \newcommand{\switchargs}[2]{#2#1}

\noindent\tt
eslahchi@karun.ipm.ac.ir\\
ghebleh@karun.ipm.ac.ir\\
hhaji@karun.ipm.ac.ir
\end{document}